\documentclass[12pt]{article}

\usepackage[margin=1in]{geometry}  % set the margins to 1in on all sides
\usepackage{graphicx}              % to include figures
\usepackage{amsmath}               % great math stuff
\usepackage{amsfonts}              % for blackboard bold, etc
\usepackage{amsthm}                % better theorem environments

\newtheorem*{thma}{Theorem A}
\newtheorem*{thmb}{Theorem B}
\newtheorem*{thmc}{Theorem C}
\newtheorem*{thmd}{Theorem D}
\newtheorem*{thme}{Theorem E}

\newtheorem*{thm1}{Theorem 1}
\newtheorem*{thm2}{Theorem 2}
\newtheorem*{thm3}{Theorem 3}
\newtheorem*{cor1}{Corollary 1}

\theoremstyle{definition}

\newtheorem*{exam1}{Example 1}
\newtheorem*{exam2}{Example 2}

\begin{document}

%\pagestyle{myheadings}
%\markboth{Eze R Nwaeze}{On the number of zeros of a Polynomial in a specified disk}
\nocite{*}

\title{\bf Some Generalizations of the Enestr\"om-Kakeya Theorem}
\date{}
\author{Eze R. Nwaeze\\
Department of Mathematics and Statistics \\
Auburn University \\
Auburn, AL 36849, USA\\
e-mail: ern0002@auburn.edu}

\maketitle

\begin{abstract}
Let $p(z)=a_0+a_1z+a_2z^2+a_3z^3+\cdots+a_nz^n$  be a polynomial of degree $n,$ where the coefficients $a_j,$ $j \in \{0,1,2,\cdots n\},$ are real numbers.  We impose some restriction on the coefficients and then prove some extensions and generalizations of the Enestr\"{o}m-Kakeya Theorem.
\end{abstract}


\section{Introduction}
A classical result due to Enestr\"om \cite{a1} and Kakeya \cite{a2} concerning the bounds for the moduli of zeros of polynomials having positive coefficients is often stated as:
\begin{thma}\label{thma}
Let $p(z)=\displaystyle{\sum_{j=0}^{n}a_jz^j}$ be a polynomial with real coefficients satisfying $$0<a_0\leq a_1\leq a_2\leq a_3\ldots\leq a_n.$$ Then all the zeros of $p(z)$ lie in $|z|\leq 1.$
\end{thma}
In the literature there exist several extensions and generalizations of this result (see \cite{a3}, \cite{a4}, \cite{a5} and \cite{a6}). Joyal et al. \cite{a5} extended Theorem A to the polynomials whose coefficients are monotonic but not necessarily nonnegative. In fact, they proved the following result:

\begin{thmb}\label{thm1.2}
Let $p(z)=\displaystyle{\sum_{j=0}^{n}a_jz^j}$ be a polynomial degree $n,$ with real coefficients satisfying $$a_0\leq a_1\leq a_2\leq a_3\ldots\leq a_n.$$ Then all the zeros of $p(z)$ lie in the disk $$|z|\leq \dfrac{1}{|a_n|}(a_n-a_0+|a_0|).$$
\end{thmb}

Aziz and Zargar \cite{a8} relaxed the hypothesis in several ways and among other things proved the following result:

\begin{thmc}
Let $p(z)=\displaystyle{\sum_{j=0}^{n}a_jz^j}$ be a polynomial of degree $n$ such that for some $k\geq1,$ $$0<a_0\leq a_1\leq a_2\leq a_3\ldots\leq ka_n.$$  Then all the zeros of $p(z)$ lie in the disk $$|z+k-1|\leq k.$$
\end{thmc}

In 2012, they further generalized Theorem C which is  an interesting extension of Theorem A. In particular, they \cite{a7} proved the following results:
\begin{thmd}
Let $p(z)=\displaystyle{\sum_{j=0}^{n}a_jz^j}$ be a polynomial of degree $n.$ If for some positive numbers $k$ and $\rho$ with  $k\geq1,$  $0<\rho\leq 1, $  $$0\leq \rho a_0\leq a_1\leq a_2\leq a_3\ldots\leq ka_n,$$  then all the zeros of $p(z)$ lie in the disk $$|z+k-1|\leq k+\dfrac{2a_0}{a_n}(1-\rho).$$

\end{thmd}

\begin{thme}
Let $p(z)=\displaystyle{\sum_{j=0}^{n}a_jz^j}$ be a polynomial of degree $n.$ If for some positive number $\rho,$  $0<\rho\leq 1, $ and for some nonnegative integer $\lambda,$  $ 0\leq \lambda< n,$~$$\rho a_0\leq a_1\leq a_2\leq a_3\ldots\leq a_{\lambda-1}\leq a_{\lambda}\geq a_{\lambda+1}\geq\cdots\geq a_{n-1}\geq a_n,$$
then all the zeros of $p(z)$ lie in the disk $$\Big|z+\dfrac{a_{n-1}}{a_n}-1\Big|\leq \dfrac{1}{|a_n|}\Big[2a_{\lambda}-a_{n-1}+(2-\rho)|a_0|-\rho a_0\Big].$$.
\end{thme}

Looking at Theorem D, one might want to know what happens if $\rho a_{0}$ is NOT nonnegative. In this paper we prove some extensions and generalization of Theorems D and E which in turns gives an answer to our enquiry.

\section{Main Results}
\begin{thm1}
Let $p(z)=\displaystyle{\sum_{j=0}^{n}a_jz^j}$ be a polynomial of degree $n.$ If for some real numbers $\alpha$ and $\beta,$ $$a_0-\beta\leq a_1\leq a_2\leq\ldots\leq a_n+\alpha,$$  then all the zeros of $p(z)$ lie in the disk $$\Big|z+\dfrac{\alpha}{a_n}\Big|\leq \dfrac{1}{|a_n|}\Big[a_n+\alpha-a_0+\beta+|\beta|+|a_0|\Big].$$
\end{thm1}

 If  $\alpha=(k-1)a_n$ ~with $k\geq 1$ and $\beta=(1-\rho)a_0$~ with $0<\rho\leq 1,$ we get the following.
\begin{cor1}
Let $p(z)=\displaystyle{\sum_{j=0}^{n}a_jz^j}$ be a polynomial of degree $n.$ If for some postive numbers $k\geq 1$ and $\rho,$ with $0<\rho\leq 1,$  $$\rho a_0\leq a_1\leq a_2\leq\ldots\leq k a_n,$$  then all the zeros of $p(z)$ lie in the disk $$|z+k-1|\leq \dfrac{1}{|a_n|}\Big[(ka_n-\rho a_0)+|a_0|(2-\rho)\Big].$$
\end{cor1}

If $a_0>0,$ then Corollary 1 amounts to Theorem D.

%\begin{cor2}\label{cor2}
%Let $p(z)=\displaystyle{\sum_{j=0}^{n}a_jz^j}$ be a polynomial of degree $n.$ If for some postive numbers $\alpha,$ ~$\rho,$ with $0<\rho\leq 1,$ $$0<\rho a_0\leq a_1\leq a_2\leq\ldots\leq a_n+\alpha.$$  Then all the zeros of $p(z)$ lie in the disk $$\Big|z+\dfrac{\alpha}{a_n}\Big|\leq \dfrac{1}{|a_n|}\Big[a_n+\alpha-a_0+2\beta+|a_0|\Big].$$
%
%\end{cor2}
%
%If, also, $\lambda =n,$ then Corollary \ref{cor1} amounts to.
%
%\begin{cor3}\label{cor3}
%Let $p(z)=\displaystyle{\sum_{j=0}^{n}a_jz^j},$ $a_0\neq 0,$ be a polynomial of degree $n$ with complex coefficients. If $Re(a_j)=\alpha_j$ and $Im(a_j)=\beta_j$ for all $j.$ Suppose $$\alpha_{n}\geq\alpha_{n - 1}\geq\alpha_{n - 2}\geq \ldots\geq \alpha_{1}\geq \alpha_{0},$$
%then the number of zeros of $p$ in $|z|\leq\delta,$ $0<\delta<1,$ does not exceed $$\dfrac{1}{\log 1/\delta}\log\dfrac{M_1}{|a_0|}$$
%where $$M_1=|\alpha_0|-\alpha_0+|\alpha_n|+\alpha_n+2\displaystyle{\sum_{j=0}^{n}|\beta_j|}.$$
%
%\end{cor3}

\begin{thm2}
Let $p(z)=\displaystyle{\sum_{j=0}^{n}a_jz^j}$ be a polynomial of degree $n.$ If for some real number $s$ and for some positive integer $\lambda,$  $ 0<\lambda< n$~$$a_0-s\leq a_1\leq a_2\leq a_3\ldots\leq a_{\lambda-1}\leq a_{\lambda}\geq a_{\lambda+1}\geq\cdots\geq a_{n-1}\geq a_n,$$
then all the zeros of $p(z)$ lie in the disk $$\Big|z+\dfrac{a_{n-1}}{a_n}-1\Big|\leq \dfrac{1}{|a_n|}\Big[2a_{\lambda}-a_{n-1}+s-a_0+|s|+|a_0|\Big].$$.
\end{thm2}

If we take $s=(1-\rho)a_0,$ with $0<\rho\leq 1,$ then Theorem 2 becomes Theorem E. Instead of proving Theorem 2, we shall prove a more general case. In fact, we prove the following result:

\begin{thm3}
Let $p(z)=\displaystyle{\sum_{j=0}^{n}a_jz^j}$ be a polynomial of degree $n.$ If for some real number $t,$ $s$ and for some positive integer $\lambda,$  $ 0< \lambda<n$~$$a_0-s\leq a_1\leq a_2\leq a_3\ldots\leq a_{\lambda-1}\leq a_{\lambda}\geq a_{\lambda+1}\geq\cdots\geq a_{n-1}\geq a_n+t,$$
then all the zeros of $p(z)$ lie in the disk $$\Big|z+\dfrac{a_{n-1}}{a_n}-\Big(1+\dfrac{t}{a_n}\Big)\Big|\leq \dfrac{1}{|a_n|}\Big[2a_{\lambda}-a_{n-1}+s-a_0+|s|+|a_0|+|t|\Big].$$.
\end{thm3}

%\section{Lemma}
%For the proof of our result we shall make use of the following result (see page 171 of the second edition) \cite{Tit}.
%
%\begin{lem}
%Let $F(z)$ be analytic in $|z|\leq R.$ Let $|F(z)|\leq M$ in the disk $|z|\leq R$ and suppose $F(0)\neq 0.$ Then for $0<\delta<1$ the number of zeros of $F(z)$ in the disk $|z|\leq \delta R$ is less than $$\dfrac{1}{\log 1/\delta}\log\dfrac{M}{|F(0)|}.$$
%\end{lem}

\section{Proof of the Theorems}
\begin{proof}[Proof of Theorem 1]
Consider the polynomial
\begin{align*}
g(z)&=(1-z)p(z)\\\\
&=-a_{n}z^{n+1}+(a_n-a_{n-1})z^n+(a_{n-1}-a_{n-2})z^{n-1}+\cdots+(a_1-a_{0})z+a_0\\\\
&=-a_{n}z^{n+1}-\alpha z^n+(a_n+\alpha-a_{n-1})z^n+(a_{n-1}-a_{n-2})z^{n-1}+\cdots+(a_1-a_{0}+\beta)z-\beta z+a_0\\\\
&=-z^n(a_{n}z+\alpha)+(a_n+\alpha-a_{n-1})z^n+(a_{n-1}-a_{n-2})z^{n-1}+\cdots+(a_1-a_{0}+\beta)z-\beta z+a_0\\\\
&=-z^n(a_{n}z+\alpha)+\phi(z),
\end{align*}

where $$\phi(z)=(a_n+\alpha-a_{n-1})z^n+(a_{n-1}-a_{n-2})z^{n-1}+\cdots+(a_1-a_{0}+\beta)z-\beta z+a_0.$$
Now for $|z|=1,$ we have
\begin{align*}
|\phi(z)|& \leq |a_n+\alpha-a_{n-1}|+|a_{n-1}-a_{n-2}|+\cdots+|a_1-a_{0}+\beta|+|\beta|+|a_0|\\\\
&=a_n+\alpha-a_{n-1}+a_{n-1}-a_{n-2}+\cdots+a_1-a_{0}+\beta+|\beta|+|a_0| \\\\
&=a_n+\alpha-a_{0}+\beta+|\beta|+|a_0|.
\end{align*}
Since this is true for all complex numbers with a unit modulus, then it must also be true for $1/z.$ With this in mind, we have
\begin{equation}
|z^n\phi(1/z)|\leq a_n+\alpha-a_{0}+\beta+|\beta|+|a_0|  ~~~~~\forall z: |z|=1.
\end{equation}
Also, the function $\Phi(z)=z^n\phi(1/z)$ is analytic in $|z|\leq 1,$ hence, Inequality (1) holds inside the unit circle by the Maximum Modulus Theorem. That is,
$$|\phi(1/z)|\leq \dfrac{a_n+\alpha-a_{0}+\beta+|\beta|+|a_0|}{|z|^n} ~~~~~\forall z: |z|\leq1.$$
Replacing $z$ by $1/z,$ we get
$$|\phi(z)|\leq \Big[a_n+\alpha-a_{0}+\beta+|\beta|+|a_0|\Big]|z|^n ~~~~~\forall z: |z|\geq1.$$
Now for $|z|\geq 1,$ we obtain
\begin{align*}
|g(z)|&=|-z^n(a_{n}z+\alpha)+\phi(z)|\\\\
&\geq |z^n||a_{n}z+\alpha|-|\phi(z)|\\\\
&\geq  |z^n||a_{n}z+\alpha|- \Big[a_n+\alpha-a_{0}+\beta+|\beta|+|a_0|\Big]|z|^n\\\\
&= |z^n|\Big(|a_{n}z+\alpha|- \Big[a_n+\alpha-a_{0}+\beta+|\beta|+|a_0|\Big]\Big)\\\\
&>0
\end{align*}
if and only if $$|a_{n}z+\alpha|>\Big[a_n+\alpha-a_{0}+\beta+|\beta|+|a_0|\Big]$$
if and only if $$\Big|z+\dfrac{\alpha}{a_n}\Big|>\dfrac{1}{|a_n|}\Big[a_n+\alpha-a_{0}+\beta+|\beta|+|a_0|\Big].$$
Thus, all the zeros of $g(z)$ whose modulus is greater than or equal to 1 lie in
\begin{equation}\label{eqn2}
\Big|z+\dfrac{\alpha}{a_n}\Big|\leq\dfrac{1}{|a_n|}\Big[a_n+\alpha-a_{0}+\beta+|\beta|+|a_0|\Big].
\end{equation}

But those zeros of $p(z)$ whose modulus is less than 1 already satisfy (\ref{eqn2}) - since $|\phi(z)|\leq a_n+\alpha-a_{0}+\beta+|\beta|+|a_0|  ~~~~~\forall z: |z|= 1$ and $\phi(z)=g(z)+z^n(a_{n}z+\alpha)$. Also, all the zeros of $p(z)$ are zeros of $g(z).$ That completes the proof of Theorem 1.
\end{proof}

\begin{proof}[Proof of Theorem 3]
Consider the polynomial
\begin{align*}
g(z)&=(1-z)p(z)\\\\
&=-a_{n}z^{n+1}+(a_n-a_{n-1})z^n+(a_{n-1}-a_{n-2})z^{n-1}+\cdots+(a_1-a_{0})z+a_0\\\\
&=-a_{n}z^{n+1}+(a_n-a_{n-1})z^n+(a_{n-1}-a_{n-2})z^{n-1}+\cdots+(a_{\lambda+1}-a_{\lambda})z^{\lambda+1}\\\\
&+(a_{\lambda}-a_{\lambda-1})z^{\lambda}+\cdots+ (a_1-a_{0})z+a_0\\\\
&=-z^n[a_{n}z-a_n+a_{n-1}-t]-tz^n+(a_{n-1}-a_{n-2})z^{n-1}+\cdots+(a_{\lambda+1}-a_{\lambda})z^{\lambda+1}\\\\
&+(a_{\lambda}-a_{\lambda-1})z^{\lambda}+\cdots+ (a_1-a_{0}+s)z-sz+a_0\\\\
&=-z^n[a_{n}z-a_n+a_{n-1}-t]+\psi(z),
\end{align*}
where $$\psi(z)=-tz^n+(a_{n-1}-a_{n-2})z^{n-1}+\cdots+(a_{\lambda+1}-a_{\lambda})z^{\lambda+1}
+(a_{\lambda}-a_{\lambda-1})z^{\lambda}+\cdots+ (a_1-a_{0}+s)z-sz+a_0.$$

For $|z|=1,$ we get
\begin{align*}
|\psi(z)|&\leq |t|+|a_{n-1}-a_{n-2}|+\cdots+|a_{\lambda+1}-a_{\lambda}|
+|a_{\lambda}-a_{\lambda-1}|+\cdots+ |a_1-a_{0}+s|+|s|+|a_0|\\\\
&=|t|+a_{n-2}-a_{n-1}+\cdots+a_{\lambda}-a_{\lambda+1}
+a_{\lambda}-a_{\lambda-1}+\cdots+ a_1-a_{0}+s+|s|+|a_0|\\\\
&=|t|-a_{n-1}+2a_{\lambda}-a_{0}+s+|s|+|a_0|.
\end{align*}
It is clear that
\begin{equation}
|z^n\psi(1/z)|\leq|t|-a_{n-1}+2a_{\lambda}-a_{0}+s+|s|+|a_0|
\end{equation}
on the unit circle. Since the function $\Psi(z)=z^n\psi(1/z)$ is analytic in $|z|\leq 1,$ Inequality (3)  holds inside the unit circle by the Maximum Modulus Theorem. That is,
$$|\psi(1/z)|\leq\dfrac{|t|-a_{n-1}+2a_{\lambda}-a_{0}+s+|s|+|a_0|}{|z|^n} $$
for $|z|\leq 1.$ Replacing $z$ by $1/z$ we get
$$|\psi(z)|\leq\Big[|t|-a_{n-1}+2a_{\lambda}-a_{0}+s+|s|+|a_0|\Big]|z|^n $$ for $|z|\geq 1.$\\
Now for $|z|\geq 1,$ we have
\begin{align*}
|g(z)|&\geq |z^n||a_{n}z-a_n+a_{n-1}-t|-|\psi(z)|\\\\
&\geq  |z^n||a_{n}z-a_n+a_{n-1}-t|-\Big[|t|-a_{n-1}+2a_{\lambda}-a_{0}+s+|s|+|a_0|\Big]|z|^n\\\\
&= |z^n|\Big(|a_{n}z-a_n+a_{n-1}-t|-\Big[|t|-a_{n-1}+2a_{\lambda}-a_{0}+s+|s|+|a_0|\Big]\Big)\\\\
&>0
\end{align*}

if and only if
$$|a_{n}z-a_n+a_{n-1}-t|>\Big[|t|-a_{n-1}+2a_{\lambda}-a_{0}+s+|s|+|a_0|\Big]$$
if and only if $$\Big|z+\dfrac{a_{n-1}}{a_n}-\Big(1+\dfrac{t}{a_n}\Big)\Big|>\dfrac{1}{|a_n|}\Big[|t|-a_{n-1}+2a_{\lambda}-a_{0}+s+|s|+|a_0|\Big].$$

Hence, the zeros of $p(z)$ with modulus greater or equal to 1 are in the closed disk $$\Big|z+\dfrac{a_{n-1}}{a_n}-\Big(1+\dfrac{t}{a_n}\Big)\Big|\leq\dfrac{1}{|a_n|}\Big[|t|-a_{n-1}+2a_{\lambda}-a_{0}+s+|s|+|a_0|\Big].$$ Also, those zeros of $p(z)$ whose modulus is less than 1 already satisfy the above inequality since $\psi(z)=g(z)+z^n[a_{n}z-a_n+a_{n-1}-t]$,~ and~ for ~$|z|=1,$ $|\psi(z)|\leq |t|-a_{n-1}+2a_{\lambda}-a_{0}+s+|s|+|a_0|.$  That completes the proof.
\end{proof}

\section{Demonstrating Examples}
\begin{exam1}
Let's consider the polynomial $$p(z)=3z^5+4z^4+3z^3+2z^2+z-1.$$ The coefficients here are $a_5=3,$ $a_4=4,$ $a_3=3,$ $a_2=2,$ $a_1=1$ and $a_0=-1.$ We cannot apply Theorems A, B, C and  D. But we can apply Theorem 1 to determine where all  the zeros of the polynomial lie. Using MATLAB, we obtain the following zeros : $-0.9154 + 0.4962i, ~-0.9154 - 0.4962i,  ~0.0530 + 0.8845i, ~0.0530 - 0.8845i,~ 0.3916.$ Taking $\alpha=2$ and $\beta=0,$ Theorem 1 gives that all the zeros of the polynomial lie in the closed disk $|3z+2|\leq 7.$
\end{exam1}

\begin{exam2}
Next, consider $$q(z)=-z^6+2z^5+2z^4+3z^3+z^2-2.$$ The coefficients of $q(z)$ are $a_6=-1,$ $a_5=2,$ $a_4=2,$ $a_3=3,$ $a_2=1,$  $a_1=0$ and $a_0=-2.$ Using MATLAB, we obtain the following zeros: $3.0197, ~-0.7682 + 0.5814i, ~-0.7682 - 0.5814i, ~-0.0803 + 1.0233i, ~-0.0803 - 1.0233i, ~0.6773. $ Taking $\lambda=3,$ $t=1$ and $s=0,$ Theorem 3 gives that the zeros lie in $|z-2|\leq 9.$

\end{exam2}

\section{Acknowledgment}
The author is greatly indebted to the referee for his/her several useful suggestions and valuable comments.

\end{document}